\documentclass[preprint,11pt]{article}

\usepackage{graphicx,epsf,amsbsy,amsmath,amssymb}
\usepackage[latin1]{inputenc}

\font\bbfnt=msbm10

\def\bbR{\mbox{\bbfnt R}}

\newcounter{nuthm}
\newenvironment{thm}{\refstepcounter{nuthm}
    \begin{trivlist}\item[\hskip\labelsep{\bf Theorem~\thenuthm.}]\sl}
    {\rm\end{trivlist}}

\newcounter{nuprop}
\newenvironment{prop}{\refstepcounter{nuprop}
    \begin{trivlist}\item[\hskip\labelsep{\bf Proposition~\thenuprop.}]\sl}
    {\rm\end{trivlist}}

\begin{document}

\title{Invariant manifolds of complex systems}

\author{Jean-Marc Ginoux and Bruno Rossetto, \\ P.R.O.T.E.E. Laboratory, I.U.T. de Toulon,\\ Universit\'e du Sud,
B.P. 20132, 83957, La Garde cedex, France,\\ E-mail: ginoux@univ-tln.fr, rossetto@univ-tln.fr}

\date{{\bf Keywords}: Invariant curves, invariant surfaces, multiple time scales
dynamical systems, complex systems.}

\maketitle

\begin{abstract}
The aim of this work is to establish the existence of invariant
manifolds in {\it complex systems}. Considering {\it trajectory curves}
integral of multiple time scales dynamical systems of
dimension two and three (predator-prey models, neuronal
bursting models) it is shown that there exists in the phase
space a {\it curve} (resp. a {\it surface}) which is invariant with
respect to the flow of such systems. These invariant
manifolds are playing a very important role in the stability of
complex systems in the sense that they are "restoring" the
determinism of {\it trajectory curves}.
\end{abstract}

\section{Dynamical systems}

In the following we consider a system of ordinary
differential equations defined in a compact E included in
\begin{equation}
  \label{eqgi1}
  \frac{{\rm d}\vec{X}}{{\rm d}t} = \vec{\Im} \left( \displaystyle \vec{X} \right)
\end{equation}
with $\vec{X} = \left[ x_1,x_2, ...,x_n \right]^t \in E \subset \bbR^n$ and
\[
  \vec{\Im} \left( \displaystyle \vec{X} \right) =
  \left[ \displaystyle f_1 \left( \displaystyle \vec{X} \right), f_2 \left( \displaystyle \vec{X} \right), ...
     f_n \left( \displaystyle \vec{X} \right)
  \right]^t \in E \subset \bbR^n \, .
\]

The vector $\vec{\Im} \left( \displaystyle \vec{X} \right)$ defines a velocity vector field in E
whose components $f_1$ which are supposed to be continuous
and infinitely differentiable with respect to all $x_t$ and $t$, i.e., are $C^\infty$ functions in E and with
values included in $\bbR$. For more details, see for example \cite{Cod55}. A solution of this system is an
{\it integral curve} $\vec{X} (t)$ tangent to $\vec{\Im}$ whose values define the states of the
dynamical system described by Equation (\ref{eqgi1}). Since none of the components $f_i$ of the velocity
vector field depends here explicitly on time, the system is said to be autonomous.

\section{Trajectory curves}

The integral of the system (\ref{eqgi1}) can be associated with the coordinates, i.e., with the position, of a
point M at the instant $t$. The total derivative of $\vec{V}(t)$ namely the instantaneous
acceleration vector $\vec{\gamma} (t)$ may be written, while using the chain rule, as:
\begin{equation}
  \label{eqgi2}
  \vec{\gamma} = \frac{{\rm d}\vec{V}}{{\rm d}t} = \frac{{\rm d}\vec{\Im}}{{\rm d} \vec{X}}
  \frac{{\rm d}\vec{X}}{{\rm d}t} = J \vec{V}
\end{equation}
where $\frac{{\rm d}\vec{\Im}}{{\rm d} \vec{X}}$ is the functional jacobian matrix $J$ of the system
(\ref{eqgi1}). Then, the {\it integral curve} defined by the vector function $\vec{X}(t)$ of the scalar
variable $t$ representing the trajectory of M can be considered as a {\it plane} or a {\it space curve}
which has local metrics properties namely {\it curvature} and {\it torsion}.

\subsection{Curvature}

The curvature, which expresses the rate of changes of the tangent to the trajectory curve of system
(\ref{eqgi1}), is defined by
\begin{equation}
  \frac{1}{\Re} = \frac{ \| \displaystyle \vec{\gamma} \wedge \vec{V} \|}{\| \vec{V} \|^3}
\end{equation}
where $\Re$ represents the {\it radius of curvature}.

\subsection{Torsion}

The {\it torsion}, which expresses the difference between the {\it trajectory curve} of system (\ref{eqgi1})
and a {\it plane curve}, is defined by:
\begin{equation}
  \frac{1}{\Im} =
    - \frac{ \displaystyle \dot{\vec{\gamma}} \cdot
      \left( \displaystyle \vec{\gamma} \wedge \vec{V} \right) }{
      \| \displaystyle \vec{\gamma} \wedge \vec{V} \|^2 }
\end{equation}
where $\Im$ represents the {\it radius of torsion}.

\section{Lie Derivative --- Darboux Invariant}

Let $\varphi$ a $C^4$ function defined in a compact E included in $\bbR$ and $\vec{X}(t)$
the integral of the dynamic system defined by (\ref{eqgi1}). The Lie derivative is defined as follows:
\begin{equation}
  L_{\vec{X}} \varphi = \vec{V} \cdot \vec{\nabla} \varphi
  = \sum_{i=1}^n{\displaystyle \frac{\partial \varphi}{\partial x_i} \dot{x}_i} =
  \frac{{\rm d} \varphi}{{\rm d}t}
\end{equation}

\begin{thm}

An {\it invariant curve} (resp. {\it surface}) is defined by
$\varphi (\vec{X}) = 0$ where $\varphi$ is a $C^1$  in an open set U and such there exists a
$C^4$ function denoted $k(\vec{X})$ and called cofactor which satisfies

\begin{equation}
  L_{\vec{X}} \phi (\vec{X}) = k (\vec{X}) \phi (\vec{X})
\end{equation}
for all $\vec{X} \in U$.
\end{thm}

Proof of this theorem may be found in \cite{Dar78}.

\begin{thm}
If $L_{\vec{X}} \varphi =0$ then $\varphi$ is first integral of the dynamical system
defined by (\ref{eqgi1}). So, $\varphi$ is first integral of the dynamical system
defined by $\{ \varphi = \alpha \}$ and where $\alpha$ is constant.

\end{thm}

Proof of this theorem may be found in \cite{Dem89}.

\section{Invariant Manifolds}

According to the previous theorems 1 and 2 the following
proposition may be established.

\begin{prop}

The location of the points where the local
curvature of the trajectory curves integral of a two
dimensional dynamical system defined by (\ref{eqgi1}) vanishes is
first integral of this system. Moreover, the invariant curve
thus defined is over flowing invariant with respect to the
dynamical system (\ref{eqgi1}).

\end{prop}

Proof of this theorem may be found in \cite{Gin06}.

\begin{prop}

The location of the points where the local
torsion of the trajectory curves integral of a three
dimensional dynamical system defined by (\ref{eqgi1}) vanishes is
first integral of this system. Moreover, the invariant surface
thus defined is over flowing invariant with respect to the
dynamical system (\ref{eqgi1}).

\end{prop}

Proof of this theorem may be found in \cite{Gin06}.

\section{Applications to Complex Systems}

According to this method it is possible to show that any
dynamical system defined by (\ref{eqgi1}) possess an invariant
manifold which is endowing stability with the trajectory
curves, restoring thus the loss determinism inherent to the
non-integrability feature of these systems. So, this method
may be also applied to any complex system such that
predator-prey models, neuronal bursting models...
But, in order to give the most simple and consistent
application, let's focus on two classical examples:

\begin{itemize}

\item the Balthazar Van der Pol model;

\item the Lorenz model.

\end{itemize}

\subsection{Van der Pol model}

The oscillator of B. Van der Pol \cite{VdP26} is a second-order
system with non-linear frictions which can be written:
\[ \ddot{x} + \alpha \left( \displaystyle x^2 -1 \right) \dot{x} + x = 0 \, .
\]
The particular form of the friction which can be carried out
by an electric circuit causes a decrease of the amplitude of
the great oscillations and an increase of the small. There are
various manners of writing the previous equation like a first
order system. One of them is:
\[
  \left\{
    \begin{array}{l}
      \dot{x} = \alpha \left( \displaystyle x + y - \frac{x^3}{3} \right) \\[0.3cm]
      \dot{y} = - \displaystyle \frac{x}{\alpha}
    \end{array}
  \right.
\]
When $\alpha$ becomes very large, $x$ becomes a fast variable
and $y$ a slow variable. In order to analyze the limit $\alpha \rightarrow \infty$, we
introduce a small parameter $\epsilon = \frac{1}{\alpha^2}$ and a slow time
$t'=\frac{t}{\alpha} \sqrt{\epsilon t}$. Thus, the system can be written:
\begin{equation}
  \label{eqgi7}
  \vec{V} =
  \left( \displaystyle
    \begin{array}{l}
      \displaystyle \frac{{\rm d}x}{{\rm d}t} \\[0.3cm]
      \displaystyle \frac{{\rm d}y}{{\rm d}t}
    \end{array}
  \right) = \Im
  \left(
    \begin{array}{l}
      \displaystyle f (x,y) \\[0.3cm]
      \displaystyle g (x,y)
    \end{array}
  \right) =
  \left(
    \begin{array}{c}
      \displaystyle \frac{1}{\epsilon} \left( \displaystyle x + y - \frac{x^3}{3} \right) \\[0.3cm]
      \displaystyle -x
    \end{array}
  \right)
\end{equation}
with $\epsilon$ a positive real parameter: $\epsilon=0.05$ and where the functions
$f$ and $g$ are infinitely differentiable with respect to all
$x_i$ and $t$, i.e. are $C^\infty$ functions in a compact E included in
$\bbR^2$ and with values in $\bbR$.

According to Proposition 1, the location of the points where
the local {\it curvature} vanishes leads to the following equation:
\begin{equation}
  \label{eqgi8}
  \phi (x,y) = 9 y^2 + \left( \displaystyle 9x + 3x^3 \right) y + 6x^4 -2x^6 + 9x^2 \epsilon
\end{equation}

According to Theorem 1 (Cf. Appendix for details), the Lie
derivative of Equation (\ref{eqgi8}) may be written:
\begin{equation}
  \label{eqgi9}
  L_{\vec{X}} \phi (\vec{X}) = \mbox{Tr} [J] \phi (\vec{X}) + \frac{2x^2}{\epsilon}
  \left( -3x -3y + x^3 \right)
\end{equation}

Let's plot the function $\phi (x,y)$ (in blue), its Lie derivative $L_{\vec{X}} \phi (\vec{X})$
(in magenta), the {\it singular approximation} $x+ y - \frac{x^3}{3}$ (in green) and the {\it limit cycle}
corresponding to system (\ref{eqgi7}) (in red):

\begin{figure}[ht]
\begin{center}
\includegraphics[height=10.5cm]{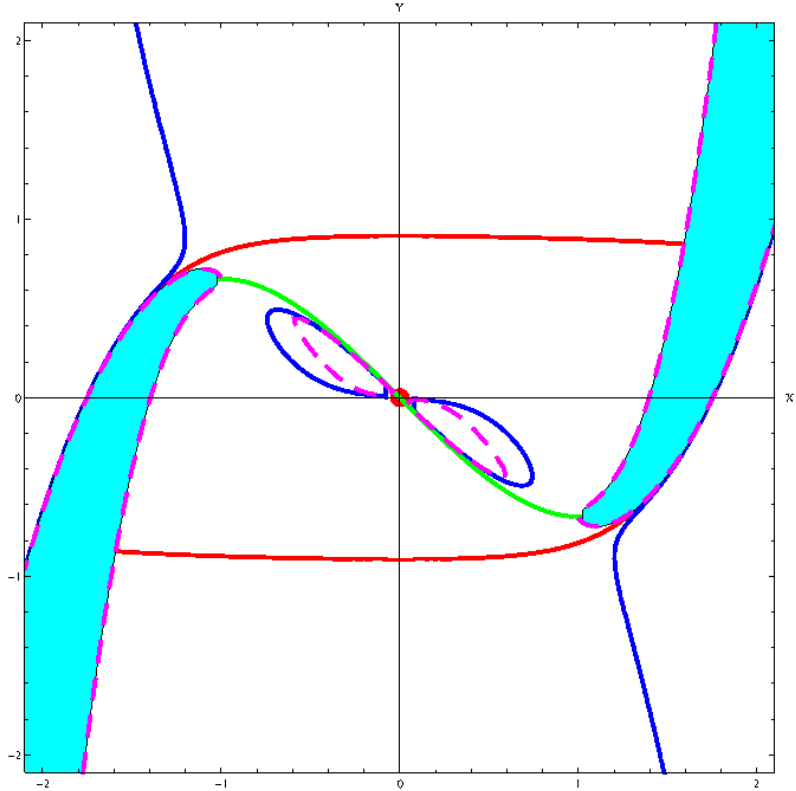}
\caption{Van der Pol model.}
\label{figin1}
\end{center}
\end{figure}

According to Fenichel's theory, there exists a function $\varphi (x,y)$ defining a manifold
(curve) which is overflowing invariant and which is $C^r {\cal O}(\epsilon)$ close to the {\it singular
approximation}. It is easy to check that in the vicinity of the {\it singular approximation} which
corresponds to the second term of the right-hand-side of Equation (\ref{eqgi9}) we have:

\[ L_{\vec{X}} \phi (\vec{X}) = \mbox{Tr}[J] \phi (\vec{X}) \, .  \]
Moreover, it can be shown that in the location of the points
where the local {\it curvature} vanishes, i.e., where $\varphi (x,y)=0$.
Equation (\ref{eqgi9}) can be written:
\[ L_{\vec{X}} \phi (\vec{X}) = 0 \, . \]
So, according to Theorem 1 and 2, we can claim that the
manifold defined by $\varphi (x,y)=0$ is an {\it invariant curve} with
respect to the flow of system (\ref{eqgi7}) and is a {\it local first integral}
of this system.

\subsection{Lorenz model}

The purpose of the model established by Edward Lorenz \cite{Lor63}
was in the beginning to analyze the impredictible behaviour of weather. It most widespread form is as
follows:
\begin{equation}
  \label{eqgi10}
  \vec{V} =
  \left( \displaystyle
    \begin{array}{l}
      \displaystyle \frac{{\rm d}x}{{\rm d}t} \\[0.3cm]
      \displaystyle \frac{{\rm d}y}{{\rm d}t} \\[0.3cm]
      \displaystyle \frac{{\rm d}z}{{\rm d}t}
    \end{array}
  \right) = \Im
  \left(
    \begin{array}{l}
      \displaystyle f (x,y,z) \\[0.3cm]
      \displaystyle g (x,y,z) \\[0.3cm]
      \displaystyle h (x,y,z)
    \end{array}
  \right) =
  \left(
    \begin{array}{c}
      \displaystyle \sigma (y-x) \\[0.3cm]
      \displaystyle -xz + rx -y \\[0.3cm]
      \displaystyle xy - \beta z
    \end{array}
  \right)
\end{equation}
with $\sigma$, $r$ and $\beta$ are real parameters: $\sigma =10$, $\beta=\frac{8}{3}$,
$r = 28$ and where the functions $f$, $g$ and $h$ are infinitely differentiable with respect to all
$x_i$, and $t$, i.e., are $C^\infty$ functions in a compact E included in $\bbR^3$ and with values in
$\bbR$. According to Proposition 1, the location of the points where the local torsion vanishes leads to an
equation which for place reasons can not be expressed. Let's name it as previously:
\begin{equation}
  \label{eqgi11}
  \varphi (x,y,z) \, .
\end{equation}
According to Theorem 1 (Cf. Appendix for details), the Lie
derivative of Equation (\ref{eqgi11}) may be written:
\begin{equation}
  \label{eqgi12}
  L_{\vec{X}} \phi (\vec{X}) = \mbox{Tr}[J] \phi (\vec{X}) + P (\vec{V} \cdot \vec{\gamma})
\end{equation}
where $P$ is a polynomial function of both vectors $\vec{V}$ and $\vec{\gamma}$.
Let's plot the function $\phi (x,y,z)$ and its Lie derivative $L_{\vec{X}} \phi (\vec{X})$
and the {\it attractor} corresponding to system (\ref{eqgi10}):

\begin{figure}[ht]
\begin{center}
\includegraphics[height=10cm]{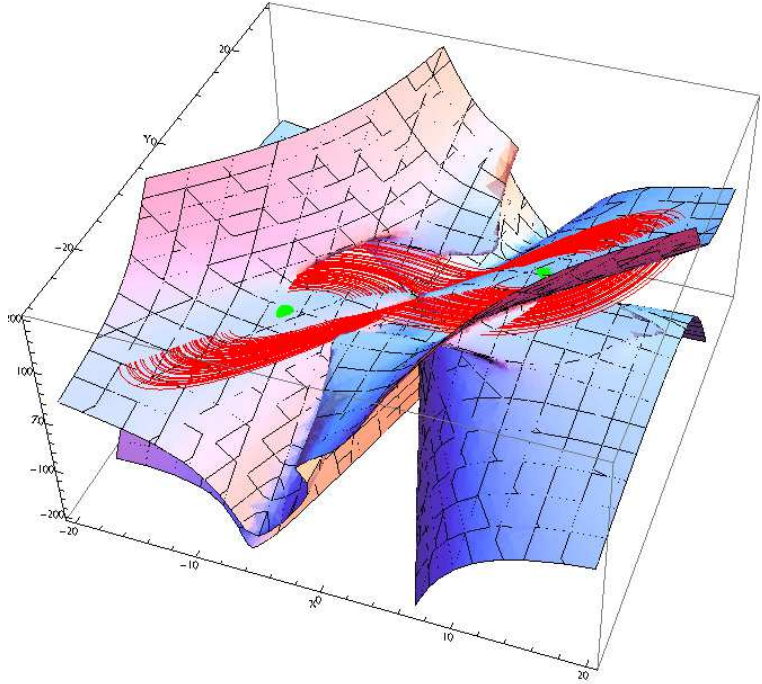}
\caption{Lorenz model.}
\label{figin2}
\end{center}
\end{figure}

It is obvious that the function $\phi (x,y,z)$ defining a manifold (surface) is merged into
the corresponding to its Lie derivative. It is easy to check that in the vicinity of the
manifold $\phi (x,y,z)$ Equation (\ref{eqgi12}) reduces to:
\[ L_{\vec{X}} \phi (\vec{X}) = \mbox{Tr}[J] \phi (\vec{X})  \, . \]
Moreover, it can be shown that in the location of the points
where the local torsion vanishes, i.e., where $\phi (x,y,z)=0$
Equation (\ref{eqgi12}) can be written:
\[ L_{\vec{X}} \phi (\vec{X}) = 0 \]
So, according to Theorem 1 and 2, we can claim that the
manifold defined by $\phi (x,y,z)=0$ is an invariant surface
with respect to the flow of system (\ref{eqgi10}) and is a local first
integral of this system.

\section{Discussion}

In this work, existence of {\it invariant manifolds} which
represent local first integrals of two (resp. three) dimensional
dynamical systems defined by (\ref{eqgi1}) has been established.
From these two characteristics it can be stated that the
former implies that such manifolds are representing the
stable part of the trajectory curves in the phase space and
from the latter that they are restoring the loss determinism
inherent to the non-integrability feature of such systems.
Moreover, while considering that dynamical systems defined
by (\ref{eqgi1}) include complex systems, it is possible to apply this
method to various models of ecology (predator-prey
models), neuroscience (neuronal bursting models), molecular
biology (enzyme kinetics models)... Research of such
invariant manifolds in coupled systems or in systems of
higher dimension (four and more) would be of great interest.

\section*{Acknowledgements}

Authors would like to thank Professors M. Aziz-Alaoui and
C. Bertelle for their useful collaboration.

\appendix

\section*{Appendix}

First of all, let's recall the following results:

\begin{equation*}
  L_{\vec{X}} \| \vec{u} \| = \frac{{\rm d} \| \vec{u} \|}{{\rm d}t} =
  \frac{\vec{u} \cdot \dot{\vec{u}}}{\| \vec{u} \|} \quad \mbox{(A-1)}
\end{equation*}

\paragraph{Two-dimensional dynamical system}\hfill \\

Let's pose: $\varphi (\vec{X}) = \| \vec{\gamma} \wedge \vec{V} \|$.
According to (A-1) the Lie derivative of this expression may be written:
\begin{equation*}
  L_{\vec{X}} \varphi (\vec{X}) =
  \frac{{\rm d} \| \vec{\gamma} \wedge \vec{V} \|}{{\rm d}t} =
  \frac{ \left( \displaystyle \vec{\gamma} \wedge \vec{V} \right) \cdot
  \frac{\rm d}{{\rm d}t}  \left( \vec{\gamma} \wedge \vec{V} \right) }{\| \vec{\gamma} \wedge \vec{V} \|} \quad \mbox{(A-2)}
\end{equation*}

\vspace{0.1in}

where $\dfrac{\rm d}{{\rm d}t}  \left( \vec{\gamma} \wedge \vec{V} \right) = \dot{\vec{\gamma}} \wedge
\vec{V}$.\\

According to Equation (\ref{eqgi2}) the Lie derivative of the
acceleration vector may be written:
\begin{equation*}
  \dot{\vec{\gamma}} = J \vec{\gamma} +
  \frac{{\rm d}J}{{\rm d}t} \vec{V} \quad \mbox{(A-3)}
\end{equation*}

it leads to:

\begin{equation*}
  \begin{array}{rl}
    \displaystyle \frac{\rm d}{{\rm d}t}  \left( \vec{\gamma} \wedge \vec{V} \right)
    = \dot{\vec{\gamma}} \wedge \vec{V} &
    = \left( \displaystyle  J \vec{\gamma} + \frac{{\rm d}J}{{\rm d}t} \vec{V} \right) \wedge \vec{V} \\[0.4cm]
    & = J \vec{\gamma} \wedge \vec{V} + \dfrac{{\rm d}J}{{\rm d}t} \vec{V} \wedge \vec{V}
  \end{array} \quad \mbox{(A-4)}
\end{equation*}

Using the following identity:\\

\[
(J \vec{a}) \wedge \vec{b} + \vec{a} \wedge (J \vec{b}) = \mbox{Tr} (J) (\vec{a} \wedge \vec{b})
\]

\vspace{0.1in}

it can be established that:\\ 

\[
J \vec{\gamma} \wedge \vec{V} = \mbox{Tr} (J) ( \vec{\gamma} \wedge \vec{V})
\]

\vspace{0.1in}

So, expression (A-2) may be written:
\begin{equation*}
  \begin{array}{rl}
    \displaystyle
    L_{\vec{X}} \varphi (\vec{X}) = &
    \displaystyle
    \frac{1}{\| \vec{\gamma} \wedge \vec{V} \|}
    \left( \displaystyle \mbox{Tr} (J) \left( \vec{\gamma} \wedge \vec{V} \right)
    \cdot \left( \vec{\gamma} \wedge \vec{V} \right) \right. \\[0.5cm]
    & \left. \displaystyle + \left( \displaystyle \frac{{\rm d}J}{{\rm d}t} \vec{V} \wedge \vec{V} \right)
    \cdot \left( \displaystyle \vec{\gamma} \wedge \vec{V} \right) \right)
  \end{array} \quad \mbox{(A-5)}
\end{equation*}

Let's note that: $\left( \vec{\gamma} \wedge \vec{V} \right) \cdot \left( \vec{\gamma} \wedge \vec{V} \right)
= \| \vec{\gamma} \wedge \vec{V} \|^2$
and that: $\vec{\beta} = \dfrac{\vec{\gamma} \wedge \vec{V}}{\| \vec{\gamma} \wedge \vec{V} \|}$.

\vspace{0.1in}

So, equation (A-5) leads to:

\begin{equation*}
  L_{\vec{X}} \varphi (\vec{X}) = \mbox{Tr} (J) \| \vec{\gamma} \wedge \vec{V} \| +
  \left( \displaystyle \frac{{\rm d}J}{{\rm d}t} \vec{V} \wedge \vec{V} \right) \cdot \vec{\beta} \quad \mbox{(A-6)}
\end{equation*}

Since vector $\frac{{\rm d}J}{{\rm d}t} \vec{V} \wedge \vec{V}$ has a unique coordinate according
to the $\vec{\beta}$-direction and since we have posed: $ \varphi (\vec{X}) = \| \vec{\gamma} \wedge \vec{V}
\|$, expression (A-6) may finally be written:
\begin{equation*}
  L_{\vec{X}} \varphi (\vec{X}) = \mbox{Tr} (J) \varphi (\vec{X}) + \left\| \displaystyle
  \frac{{\rm d}J}{{\rm d}t} \vec{V} \wedge \vec{V} \right\| \quad \mbox{(A-7)}
\end{equation*}

\paragraph{Three-dimensional dynamical system}\hfill \\

Let's pose: $\varphi (\vec{X}) = \dot{\vec{\gamma}} \cdot \left( \vec{\gamma} \wedge \vec{V} \right)$.
The Lie derivative of this expression may be written:
\begin{equation*}
  L_{\vec{X}} \varphi (\vec{X}) =
  \frac{{\rm d} \left[ \displaystyle \dot{\gamma} \cdot \left( \displaystyle {\vec{\gamma}} \wedge \vec{V}
  \right) \right]}{{\rm d}t} \quad \mbox{(A-8)}
\end{equation*}

\vspace{0.1in}

According to $\dfrac{{\rm d} }{{\rm d}t}
\left[ \displaystyle \dot{\gamma} \cdot \left( \displaystyle {\vec{\gamma}} \wedge \vec{V} \right) \right]
= \ddot{\vec{\gamma}} \cdot \left( \vec{\gamma} \wedge \vec{V} \right)$, it leads to:

\begin{equation*}
  L_{\vec{X}} \varphi (\vec{X}) =
  \frac{{\rm d} \left[ \displaystyle \dot{\gamma} \cdot \left( \displaystyle {\vec{\gamma}} \wedge \vec{V}
  \right) \right]}{{\rm d}t} = \ddot{\vec{\gamma}} \cdot \left( \vec{\gamma} \wedge \vec{V} \right) \quad \mbox{(A-9)}
\end{equation*}

\vspace{0.1in}

The Lie derivative of expression (A-3) leads to:
\[ \ddot{\vec{\gamma}} = J \dot{\vec{\gamma}} +
  2 \frac{{\rm d}J}{{\rm d}t} \vec{\gamma} + \frac{{\rm d}^2 J}{{\rm d}t^2} \vec{V}
\]

Thus, expression (A-9) reads:
\begin{equation*}
  \begin{array}{rl}
    \displaystyle L_{\vec{X}} \varphi (\vec{X}) = &
    \displaystyle  \left( J \dot{\vec{\gamma}} \right) \cdot
    \left( \vec{\gamma} \wedge \vec{V} \right) \\[0.4cm]
    & + \left( \displaystyle 2 \frac{{\rm d}J}{{\rm d}t} \vec{\gamma} +
    \frac{{\rm d}^2 J}{{\rm d}t^2} \vec{V} \right) \cdot \left( \vec{\gamma} \wedge \vec{V} \right)
  \end{array} \quad \mbox{(A-10)}
\end{equation*}

It can also be established that:
\[ \left( \displaystyle J^2 \vec{\gamma} \right) \cdot \left( \vec{\gamma} \wedge \vec{V}\right)
  = \mbox{Tr} (J) \left( J \vec{\gamma} \right) \cdot \left( \vec{\gamma} \wedge \vec{V} \right)
\]
So, since we have posed: $\varphi (\vec{X}) = \dot{\vec{\gamma}} \cdot \left( \vec{\gamma} \wedge \vec{V}
\right)$, expression (A-10) may finally be written:
\begin{equation*}
  \begin{array}{rl}
    \displaystyle
    L_{\vec{X}} \varphi (\vec{X}) = &
    \displaystyle
    \mbox{Tr} (J) \varphi (\vec{X}) + \left( - \mbox{Tr} (J) \frac{{\rm d}J}{{\rm d}t} \vec{V} \right. \\[0.4cm]
    & \left.  \displaystyle
    J \frac{{\rm d}J}{{\rm d}t} \vec{V} + 2 \frac{{\rm d}J}{{\rm d}t} \vec{\gamma}
    + \frac{{\rm d}^2J}{{\rm d}t^2} \vec{V} \right) \cdot \left( \vec{\gamma} \wedge \vec{V} \right)
  \end{array} \quad \mbox{(A-11)}
\end{equation*}


\begin{thebibliography}{08}

\bibitem{Cod55}
Coddington, E.A. \& Levinson., N., 1955.
{\it Theory of Ordinary Differential Equations},
Mac Graw Hill, New York.

\bibitem{Dar78}
Darboux, G. 1878. M\'emoire sur les \'equations diff\'erentielles
alg\'ebriques du premier ordre et du premier degr\'e.
{\it Bull. Sci. Math.} S\'er. {\bf 2} (2), 60-96, 123-143, 151-200.

\bibitem{Dem89}
Demazure, M. 1989. {\it Catastrophes et Bifurcations},
Ellipses, Paris.

\bibitem{Fen79}
Fenichel, N. 1979. Geometric singular perturbation theory
for ordinary differential equations. {\it J. Diff. Eq.} {\bf 31}, 53-98

\bibitem{Gin06}
Ginoux, J.M. and Rossetto B. 2006. Invariant manifolds of
complex systems. to appear.

\bibitem{Lor63}
Lorenz, E. N. (1963). Deterministic non-periodic flows,
{\it J. Atmos. Sc.}, {\bf 20}, 130-141.

\bibitem{VdP26}
Van der Pol, B. 1926. On 'Relaxation-Oscillations',
{\it Phil. Mag.}, {\bf 7}, Vol. 2, 978-992.

\end{thebibliography}
\end{document}